\documentclass[preprint, 3p, times]{elsarticle}

\usepackage[latin1]{inputenc}
\usepackage{epstopdf}
\usepackage{subfig}
\usepackage{amsfonts,amsmath,amssymb,amsthm}
\usepackage{graphicx}
\usepackage[dvips]{color}
\usepackage{bm}
\usepackage{algorithm}
\usepackage{algorithmic}
\usepackage{relsize}
\usepackage{accents}

\newcommand{\RR}{\mathbb{R}}

\newtheorem{definition}{Definition}[section]

\makeatletter
\def\ps@pprintTitle{%
   \let\@oddhead\@empty
   \let\@evenhead\@empty
   \def\@oddfoot{\reset@font\hfil\thepage\hfil}
   \let\@evenfoot\@oddfoot
}
\makeatother

\begin{document}
\begin{frontmatter}



\title{On the topology preservation of Gneiting's functions in image registration}


\author{Chiara Bosica}

\author{Roberto Cavoretto}

\author{Alessandra De Rossi} 

\author{Hanli Qiao}
\ead{hanli.qiao@unito.it}
\date{}

\address{Department of Mathematics \lq\lq G. Peano\rq\rq, University of Torino, via Carlo Alberto 10, I--10123 Torino, Italy}
\begin{abstract} 
The purpose of image registration is to determine a transformation such that the transformed version of the source image is similar to the target one. In this paper we focus on landmark-based image registration using radial basis functions (RBFs) transformations, in particular on the topology preservation of compactly supported radial basis functions (CSRBFs) transformations. In \cite{Cavoretto13} the performances of Gneiting's and Wu's functions are compared with the ones of other well known schemes in image registration, as thin plate spline and Wendland's functions. Several numerical experiments and real-life cases with medical images show differences in accuracy and smoothness of the considered interpolation methods, which can be explained taking into account their topology preservation properties. Here we analyze analytically and experimentally the topology preservation performances of Gneiting's functions, comparing results with the ones obtained in \cite{Yang11}, where Wendland's and Wu's functions are considered.
\end{abstract}


\end{frontmatter}




\section{Introduction} \label{sec1}

In medical image analysis, registration is a crucial step to analyze two images taken at different times or coming from different sensors or situations. The problem of image registration can intuitively formulated in the following way: given two images, called \emph{source}  and \emph{target} images, respectively, find an appropriate transformation between the two images, such that it maps the source image onto the target one. The differences between the two images can derive from different conditions, and for analyzing them we want to make images more similar each other after transformation, for an overview see e.g. \cite{Goshtasby05,Modersitzki04,Modersitzki09,Rohr01,Scherzer06,Zitova03}. In particular, in \cite{Zitova03}, Zitov\'{a} and Flusser gave an overview of image registration, and presented different methods to solve this problem. 

Some methods for image registration are based on landmarks. The \textsl{landmark-based image registration} process is based on two finite sets of landmarks, i.e. scattered data points located on images, where each landmark of the source image has to be mapped onto the corresponding landmark of the target image (see \cite{Modersitzki04,Modersitzki09,Rohr01}). The landmark-based registration problem can thus be formulated in the context of multivariate scattered data interpolation.
In landmark-based image registration, the deformed results are sensitive to the displacements of landmarks. If the displacement of one landmark is far enough from  neighborhood landmarks, the deformation  will be large and the geometrical structure will change after transformation. In such case, topology violation could thus occur.

One of the most used methods in landmark-based image registration is the radial basis functions (RBFs) method, useful to handle various geometric deformations. RBFs are functions whose values depend on the distance between points and centers \cite{Buhmann03,Fasshauer}. This important property allows the use of RBFs  in interpolation problems, such as image registration. RBFs can be classified in two groups: (i) globally supported such as thin-plate spline (TPS) and Gaussian, and (ii) compactly supported such as Gneiting's,  Wendland's and Wu's functions, see \cite{Cavoretto13,Arad94,Gneiting02,Wendland05,Wu95}, respectively.

In general, globally supported RBFs (GSRBFs) can guarantee the bending energy be small but the deformed field will occur in the whole image after transformation; conversely, compactly supported RBFs (CSRBFs) make the influence of the deformation local: around a landmark in 2D images on a circle whereas in 3D images on a sphere. In \cite{Cavoretto13,Allasia14}, authors analyzed different computational properties of GSRBFs and CSRBFs  for landmark-based image registration.

However, if GSRBFs or CSRBFs are chosen to solve the registration problem, topology should be preserved. For some RBFs the presence of a shape parameter is an important characteristic. In fact, a large shape parameter leads to flat basis functions, whereas a small shape parameter results in peaked RBFs \cite{Fasshauer}. Therefore, we can use it to control influences on the registration result \cite{Fornefett01}. 

In \cite{Yang11} evaluation of topology preservation for different CSRBFs in case of landmark-based image registration was performed. The authors compared topology preservation of different CSRBFs using two criteria: locality parameter and positivity of determinant of the corresponding Jacobian matrices for different transformations. Locality parameter is the optimal support size of different CSRBFs under topology preservation condition. In \cite{LNCS14}, we instead evaluated the performances of topology preservation for a GSRBF family such as Mat\'{e}rn functions in landmark-based image registration.
 
The family of Gneiting's functions was proposed by Gneiting in 2002 for the first time \cite{Gneiting02}. Gneiting's functions are oscillatory compactly supported functions and it is known in scattered data interpolation that they achieve good approximation results \cite {Fasshauer}. Then, they were used also in image registration, obtaining accurate results \cite {Cavoretto13,Cavoretto12}. For this reason in this paper we analize the topology preservation of Gneiting's transformations under the two criteria given above and compare the numerical results with those obtained in the  paper \cite{Yang11} using Wendland's and Wu's functions.

The paper is organized as follows. Section 2 introduces the landmark-based image registration problem. Section 3 gives definitions of two kinds of Gneiting's functions and of the associated transformations. In Section 4 we evaluate topology preservation of Gneiting's functions in the case of one-landmark. Computations and numerical results are presented. Section 5 deals the four-landmarks schematic diagram case and the landmark-based registration in brain images. Finally, in Section 6, we report the conclusions and the future work on the topic.

\section{Landmark-based image registration problem}\label{sec2}

In this paper, we only consider the 2D case. For landmark-based image registration, we define a couple of landmark sets 
$$
\mathcal{S}_N=\{\textbf{x}_j \in \RR^2,j=1,2,...,N\}
$$
and 
$$
\mathcal{T}_N=\{\textbf{t}_j \in \RR^2,j=1,2,...,N\}
$$
corresponding to the source and target images, respectively. The registration can be described as follows.

The aim is to find a proper transformation $\textbf{R}: \hspace{0.2cm} \RR^2 \rightarrow \RR^2$ between $\mathcal{S}_N$ and $\mathcal{T}_N$, such that
\begin{equation} \label{Interpolate}
\textbf{R} (\textbf{x}_j)=\textbf{t}_j, \hspace{0.3cm} j=1,2,...,N.
\end{equation}

For image registration, we can interpolate displacements to fulfill the deformation. As we mentioned in Section \ref{sec1}, the influence of deformed field is limited by CSRBFs. Here the displacements can be displayed by a CSRBF interpolant $R_{k}: \RR^2\rightarrow \RR, \hspace{0.2cm}k=1,2$, of the form

\begin{equation}\label{Displacement}
R_k(\mathbf{x}) = \sum_{j=1}^N\alpha_{jk} \Psi\left(\parallel \mathbf{x}-\mathbf{x}_j \parallel\right),
\end{equation}
where $\Psi$ stands for a CSRBF, $r=\parallel \mathbf{x}-\mathbf{x}_j \parallel$ is the Euclidean distance between $\mathbf{x}$ and $\mathbf{x}_j$, and the coefficient $\alpha_{jk}$ can be calculated by solving two linear systems. In this way we obtain the transformation ${\bf R}$. 

Following \cite{Fasshauer}, when functions $\Psi$ are strictly positive definite, the matrix is invertible since all the eigenvalues are positive. Therefore, we have a unique solution of the two linear systems. Hence in this paper, all of the CSRBFs we consider are strictly positive definite. In the following we list two examples of Wendland's ($\varphi_{3,1}$) and Wu's ($\psi_{1,2}$) functions, i.e.,
\begin{align}
\varphi_{3,1}(r) \doteq \left(1-\frac{r}{c}\right)_{+}^4\left(4\frac{r}{c}+1\right), \nonumber
\end{align}
\begin{align}
\psi_{1,2}(r) \doteq \left(1-\frac{r}{c}\right)_{+}^4\left(1+4\frac{r}{c}+3\left(\frac{r}{c}\right)^2+\frac{3}{4}\left(\frac{r}{c}\right)^3\right), \nonumber
\end{align}
where $(\cdot)_+$ is the truncated power function and $\frac{r}{c} \leq 1$, $c$ being the support function size. We remark that the larger (smaller) $c$ is, the larger (smaller) the field is. We will use these functions to compare numerical results of topology preservation obtained by Gneiting's transformations.

\section{Gneiting's functions and transformations} \label{sec3}

In this section we introduce the definitions of Gneiting's transformations. In 2002, Gneiting obtained a family of oscillating compactly supported functions \cite{Gneiting02}. Starting with a function $\varphi_m$ that is strictly positive definite and radial on $\RR^m$, for $m\geq3$, and using turning bands operator \cite{Matheron73}, we get
\begin{align}\label{LI1}
\varphi_{m-2}(r)=\varphi_m(r)+\frac{r\varphi'_m(r)}{m-2}.
\end{align}
The latter is strictly positive definite and radial on $\RR^{m-2}$ \cite{Fasshauer}. In order to obtain Gneiting's functions, we start with Wendland's functions, for example 
\begin{align*}
\varphi_{4,1}=(1-r){_+^{l+1}}[(l+1)r+1].
\end{align*} 
Using the turning bands operator, we thus obtain the functions
\begin{align}\label{LI2}
\tau_{2,l}(r)=(1-r){_+^l}\Bigg(1+lr-\frac{(l+1)(l+4)}{2}r^2 \Bigg),
\end{align}
which are strictly positive definite and radial on $\RR^2$ provided $l\geq{7/2}$. From this family we list two specific Gneiting's functions in $C^2(\RR)$, i.e.,
\begin{align}\label{LI3}
\tau_{2,7/2}(r)\doteq\left(1-\frac{r}{c}\right){_+^{7/2}}\bigg(1+\frac{7}{2}\frac{r}{c}-\frac{135}{8}\bigg(\frac{r}{c}\bigg)^2 \bigg),  
\end{align}
\begin{align}\label{LI4}
\tau_{2,5}(r)\doteq\left(1-\frac{r}{c}\right){_+^5}\bigg(1+5\frac{r}{c}-27\bigg(\frac{r}{c}\bigg)^2 \bigg).
\end{align}
Under the landmark-based image registration context we define Gneiting's transformation as follows.

\begin{definition} \label{Def1} 
Given a set of source landmark points $\mathcal{S}_N$= \{$\textbf{x}_j \in \RR^2,j=1,2,...,N$\}, and the corresponding set of target landmark points $\mathcal{T}_N$=\{$\textbf{t}_j \in \RR^2,j=1,2,...,N$\}, Gneiting's transformation $\textbf{G} :\RR^2 \rightarrow \RR^2$ is such that each its component
\begin{align*}
G_k(\mathbf{x}) : \RR^2 \rightarrow \RR,\quad  k=1,2,
\end{align*}
assumes the following form
\begin{align}
G_k(\mathbf{x}) = G_k(x_1,x_2) = \sum_{j=1}^N \alpha_{jk}\tau_{2,l}\big(\parallel \mathbf{x}-\mathbf{x}_j \parallel_2\big),  \nonumber
\end{align}
with $\textbf{x}$ = $(x_1,x_2)$ and $\textbf{x}_j$ = $(x_{j1},x_{j2})$ $\in \RR^2.$
\end{definition}

According to Definition \ref{Def1}, the transformation function $G_k(\mathbf{x}) :\RR^2 \rightarrow \RR$ is calculated for each $k=1,2,$ and the coefficients $\alpha_{jk}$ are obtained by solving two systems of linear equations.

\section{Topology preservation: the one-landmark case} \label{sec4}

From \cite{Fornefett01}, under injectivity of map, the necessary conditions to preserve topology are that  the function $\textbf{H}: \RR^2 \rightarrow \RR^2$ is continuous and the Jacobian determinant is positive at each point.

In this case, the source landmark $\textbf{p}$ is shifted by $\Delta_x$ along the $x$-axis direction and by $\Delta_y$ along the $y$-axis direction to the target landmark $\textbf{q}$. The coordinates of transformation are
\begin{center}
$H_1(\textbf{x})=x+\Delta_x\Phi(||\textbf{x}-\textbf{p}||)$,  
\end{center}  
\begin{center}
$H_2(\textbf{x})=y+\Delta_y\Phi(||\textbf{x}-\textbf{p}||)$,
\end{center}
where $\Phi$ is any CSRBF.

The positivity of the Jacobian determinant requires 
\begin{equation}
\det(J(x,y))=1+\Delta_x\frac{\partial \Phi}{\partial x}+\Delta_y\frac{\partial \Phi}{\partial y}>0,
\end{equation}
i.e.
\begin{align*}
\Delta_x\frac{\partial \Phi}{\partial x}+\Delta_y\frac{\partial \Phi}{\partial y}>-1,
\end{align*}
or, equivalently,\\
\begin{align*}
\Delta_x\frac{\partial \Phi}{\partial r}\cos\theta+\Delta_y\frac{\partial \Phi}{\partial r}\sin\theta>-1,
\end{align*}
where $\Phi$ stands for $\Phi(||\textbf{x}-\textbf{p}||)$ and $r=||\textbf{x}-\textbf{p}||$. 

If we set $\Delta=\mbox{max}(\Delta_x,\Delta_y)$, the value of $\theta$ minimizing the determinant in 2D is $\frac{\pi}{4}$; thus we get
\begin{equation}\label{WC1}
\Delta\frac{\partial\Phi}{\partial r}>-\frac{1}{\sqrt{2}}.
\end{equation}
With the condition (\ref{WC1}), one can show that all principal minors of the Jacobian are positive. It follows that the transformations defined by equation (\ref{LI1}) preserve the topology if (\ref{WC1}) holds. The minimum of $\frac{\partial\Phi}{\partial r}$ depends on the localization parameter and therefore on the support size of the parameter $c$ of Gneiting's functions.

In the next subsections we compute the minimum support size of locality parameter of Gneiting's functions \eqref{LI3} and \eqref{LI4}, satisfying \eqref{WC1}. 

\subsection{Gneiting $\tau_{2,7/2}$} \label{sec4.1}
Considering Gneiting's function (\ref{LI3}) with support size $c$, we look for the minimum value of $c$ such that \eqref{WC1} is satisfied. Since
\begin{equation}\label{Gne7_1}
\frac{\partial\Phi}{\partial r}=-\frac{99r}{16c^2}\left(1-\frac{r}{c}\right)_+^{5/2}\left(8-15\frac{r}{c}\right),
\end{equation}
the value of $r$ minimizing (\ref{Gne7_1}) is
\begin{align} \label{Gne7_2}
r=4\left(\frac{29}{270}-\frac{\sqrt{301}}{270}\right)c\approx 0.17c.
\end{align}
Thus, evaluating (\ref{Gne7_1}) at (\ref{Gne7_2}) and replacing its outcome in \eqref{WC1}, we obtain that the support size must be
\begin{align*}
c>3.60\sqrt{2}\Delta\approx 5.09\Delta.
\end{align*}

\subsection{Gneiting $\tau_{2,5}$} \label{sec4.2}
Focusing now on Gneiting's function (\ref{LI4}), whose support size is always $c$, we search for the minimum $c$ satisfying (\ref{WC1}). Here we have
\begin{align}\label{Gne5_1}
\frac{\partial\Phi}{\partial r}=-\frac{21r}{c^2}\left(1-\frac{r}{c}\right)_+^4\left(4-9\frac{r}{c}\right),
\end{align}
while the value of $r$ minimizing \eqref{Gne5_1} is given by
\begin{align} \label{Gne5_2}
r=\left(\frac{19}{54}-\frac{\sqrt{145}}{54}\right)c\approx 0.13c.
\end{align}
Evaluating (\ref{Gne5_1}) at (\ref{Gne5_2}) and substituting its result in (\ref{WC1}), we obtain
\begin{align*}
c>4.43\sqrt{2}\Delta\approx 6.26\Delta.
\end{align*}

Table \ref{tab_1} gives the locality parameters for the two Gneiting's functions (\ref{LI3}) and (\ref{LI4}), which are compared with those of Wendland's and Wu's functions (see \cite{Fornefett01,Yang11}). As we mentioned in Section \ref{sec1}, the advantage of having small locality parameter is that the influence of deformed area at each landmark turns out to be small. This property allows us to have a greater local control. From Table \ref{tab_1}, we get that the locality parameter of function $\varphi_{3,1}$ and $\psi_{1,2}$ are very similar. Moreover, the latter are smaller than those of Gneiting's functions. This means that the deformed field of Wendland's and Wu's transformations in one-landmark case is similar and its corresponding transformed area is smaller than the one of Gneiting's functions. 

\begin{table} [!htb]
\centering
\caption{Minimum support size $c$ for various CSRBFs.}
\label{tab_1}       
\begin{tabular}{llll}
\hline\noalign{\smallskip}
$\varphi_{3,1}$ & $\psi_{1,2}$ & $\tau_{2,7/2}$ & $\tau_{2,5}$\\
\noalign{\smallskip}\hline\noalign{\smallskip}
$c>2.98\Delta$ & $c>2.80\Delta$ & $c>5.09\Delta$ & $c>6.26\Delta$\\ 
\noalign{\smallskip}\hline
\end{tabular}
\end{table}

\subsection{Numerical experiments}

Let us consider a grid $[0,1]\times[0,1]$ and compare then the results obtained by distortion of the grid in the shift case of the landmark $\{(0.5,0.5)\}$ in $\{(0.6,0.7)\}$. 

In Figure 4.1 we show results assuming as a support size the minimum $c$ such that (\ref{WC1}) is satisfied, while in Figure 4.2 we take a value of $c$ which does not satisfy the topology preservation condition. In both examples, $\Delta=0.2$. 

Figure 4.1 shows that for $c$ minimum, all transformations can preserve topology well. From Figure 4.2 appears instead that, if the topology preservation is not satisfied, the transformed image is deeply misrepresented above all around the shifted point. This is especially true for Gneiting's transformations, because the chosen parameter $c=0.15$ of Gneiting's functions is much smaller than locality parameters of Wendland's and Wu's functions.

\begin{figure*} \label{1_landmark_c_si}
\begin{center}
\begin{minipage}{40mm}
\includegraphics[width=4cm]{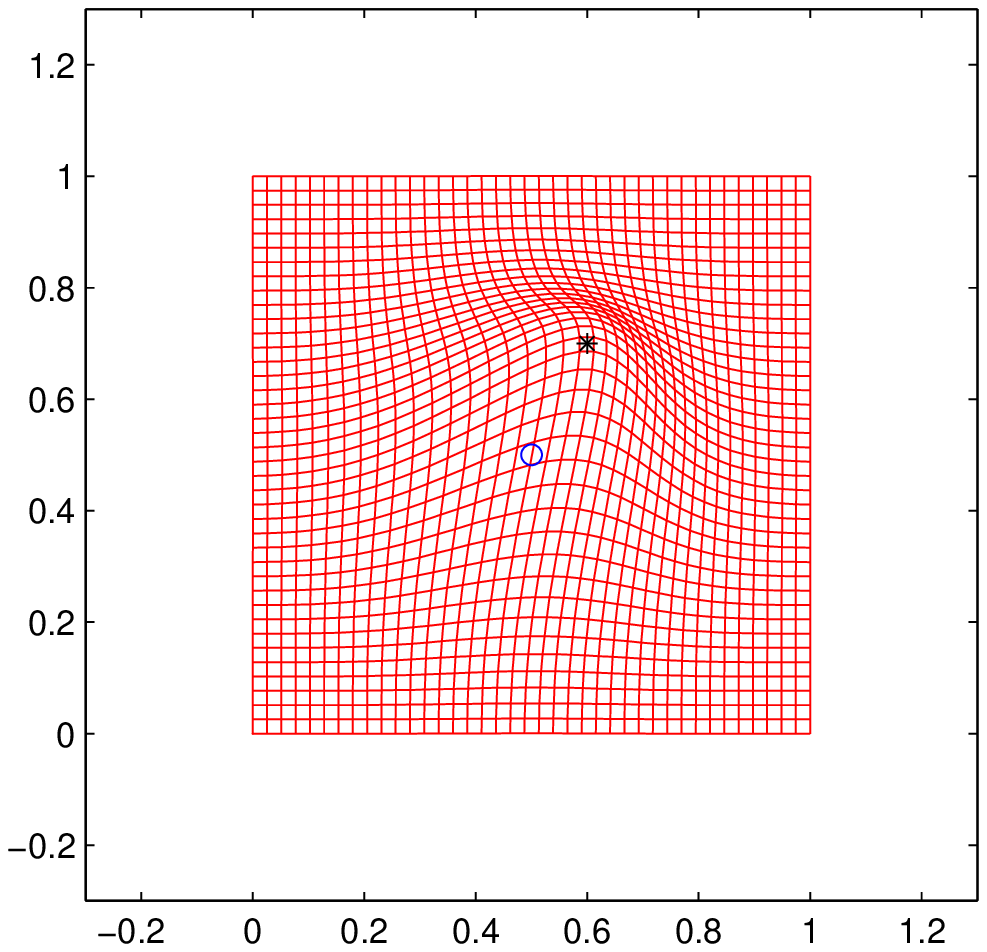}
\centerline{(a) Wendland $\varphi_{3,1}$, $c=0.6$}
\end{minipage} 
\begin{minipage}{40mm}
\includegraphics[width=4cm]{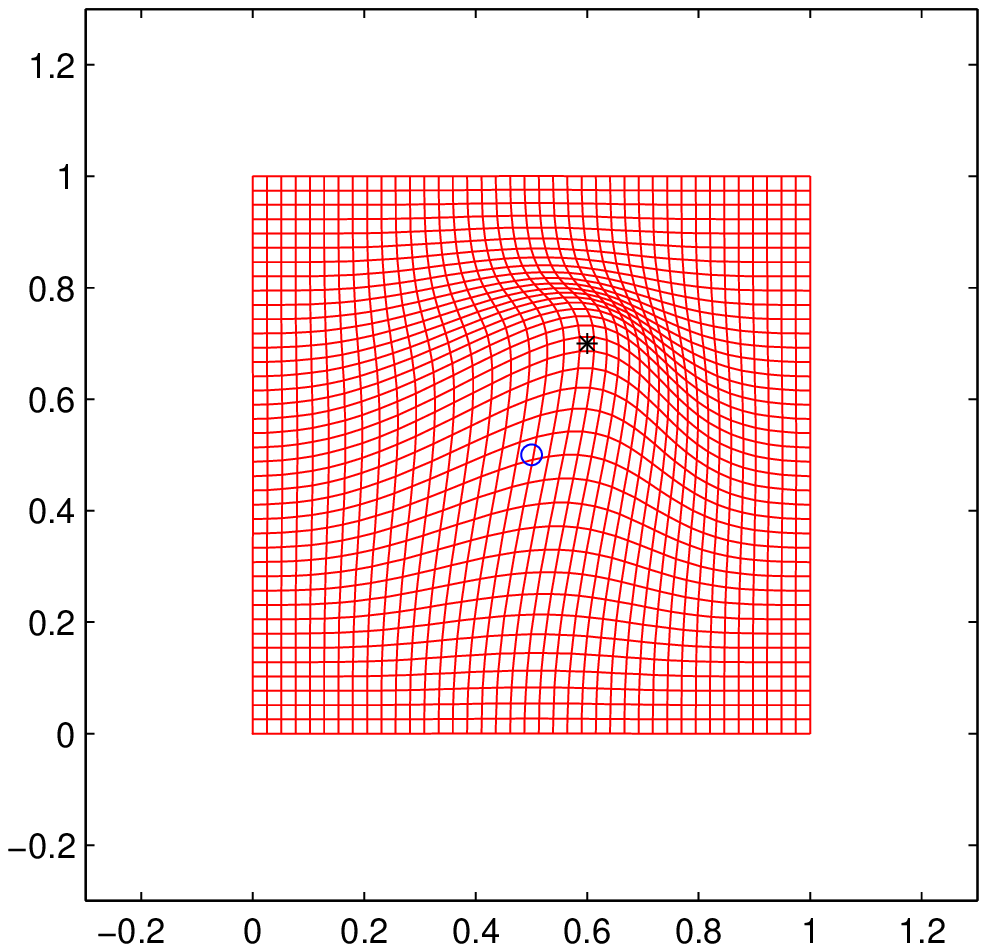}
\centerline{(b) Wu $\psi_{1,2}$, $c=0.58$}
\end{minipage} 
\begin{minipage}{40mm}
\includegraphics[width=4cm]{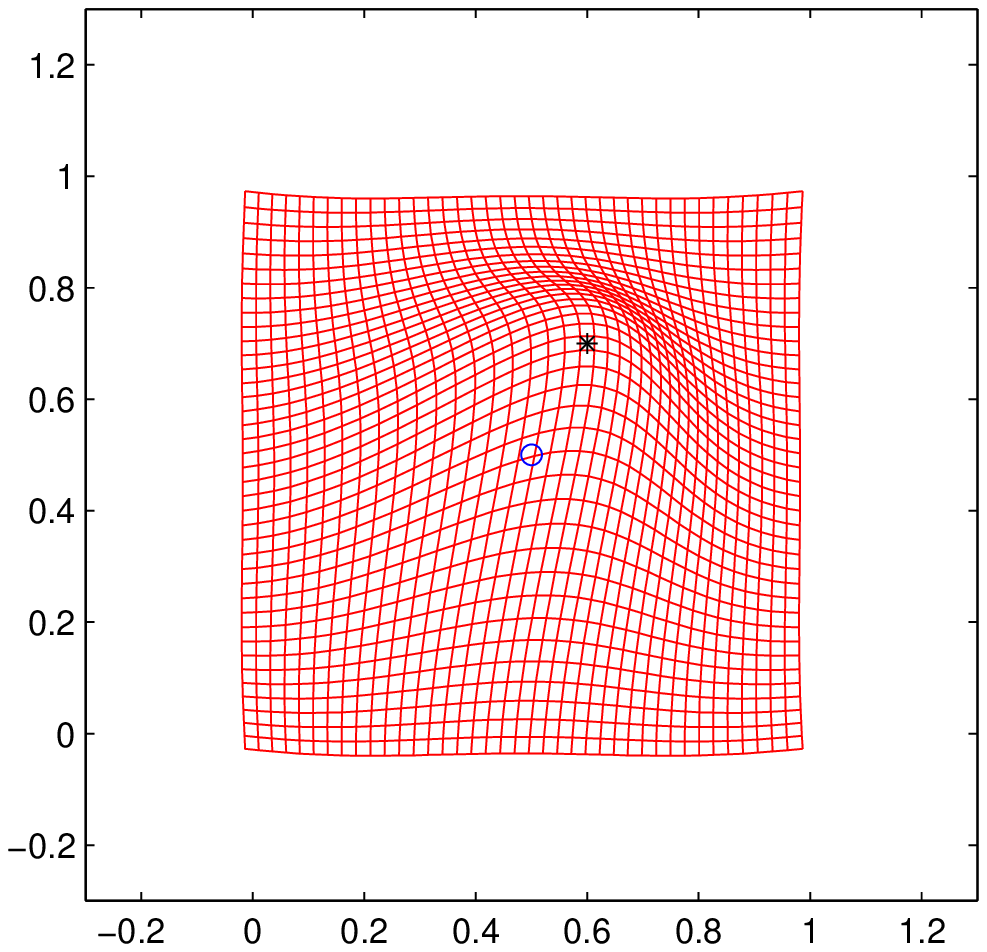}
\centerline{(c) Gneiting $\tau_{2,7/2}$, $c=1.02$}
\end{minipage} 
\begin{minipage}{40mm}
\includegraphics[width=4cm]{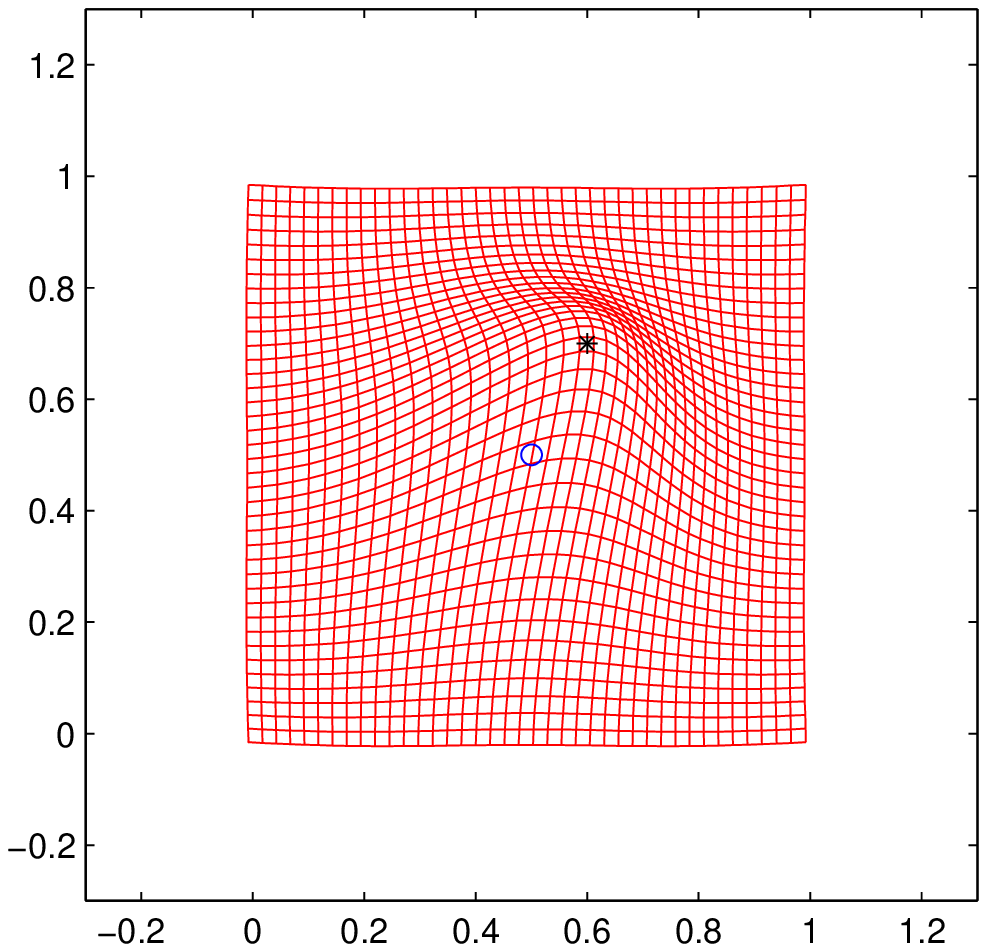}
\centerline{(d) Gneiting $\tau_{2,5}$, $c=1.26$}
\end{minipage}
\caption{Deformation results of one-landmark matching using minimum locality parameters satisfying the topology preservation condition. The source landmark is marked by a circle ($\circ$), while the target one by a star ($\ast$).}
\end{center}
\end{figure*}

\begin{figure*}
\begin{center}
\begin{minipage}{40mm}
\includegraphics[width=4cm]{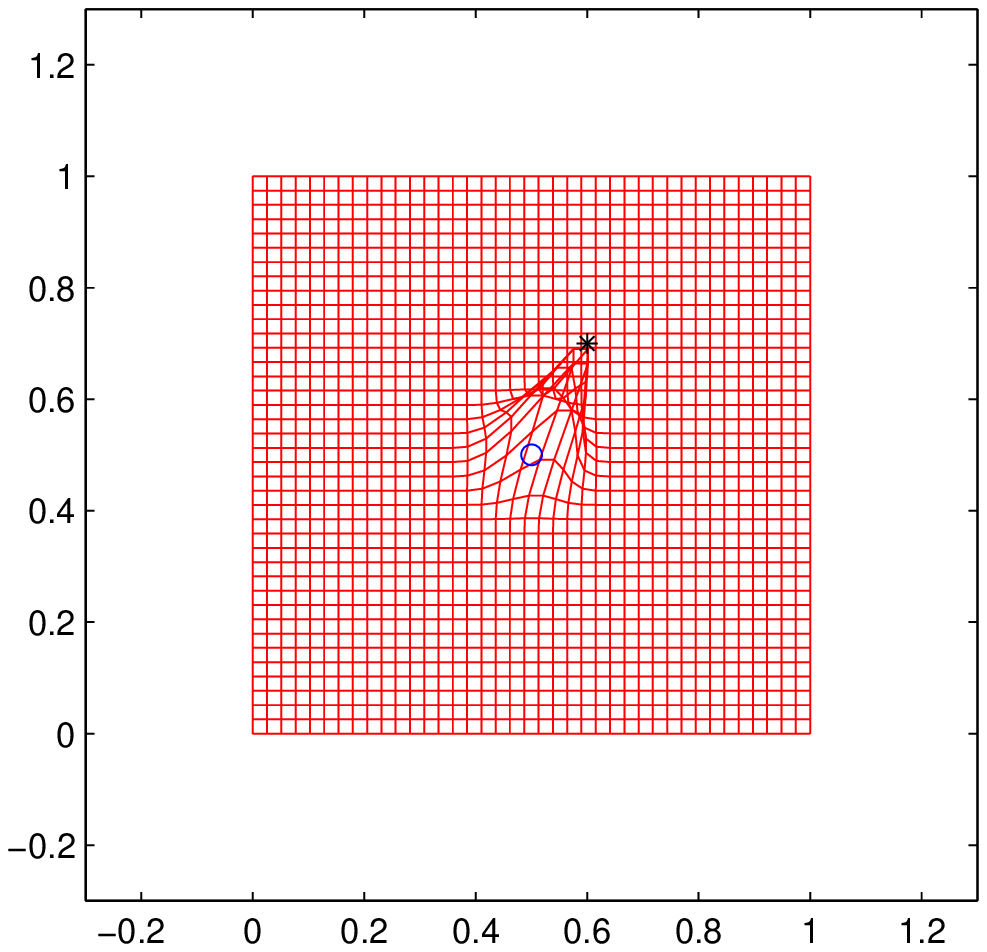}
\centerline{(a) Wendland $\varphi_{3,1}$, $c=0.15$}
\end{minipage} 
\begin{minipage}{40mm}
\includegraphics[width=4cm]{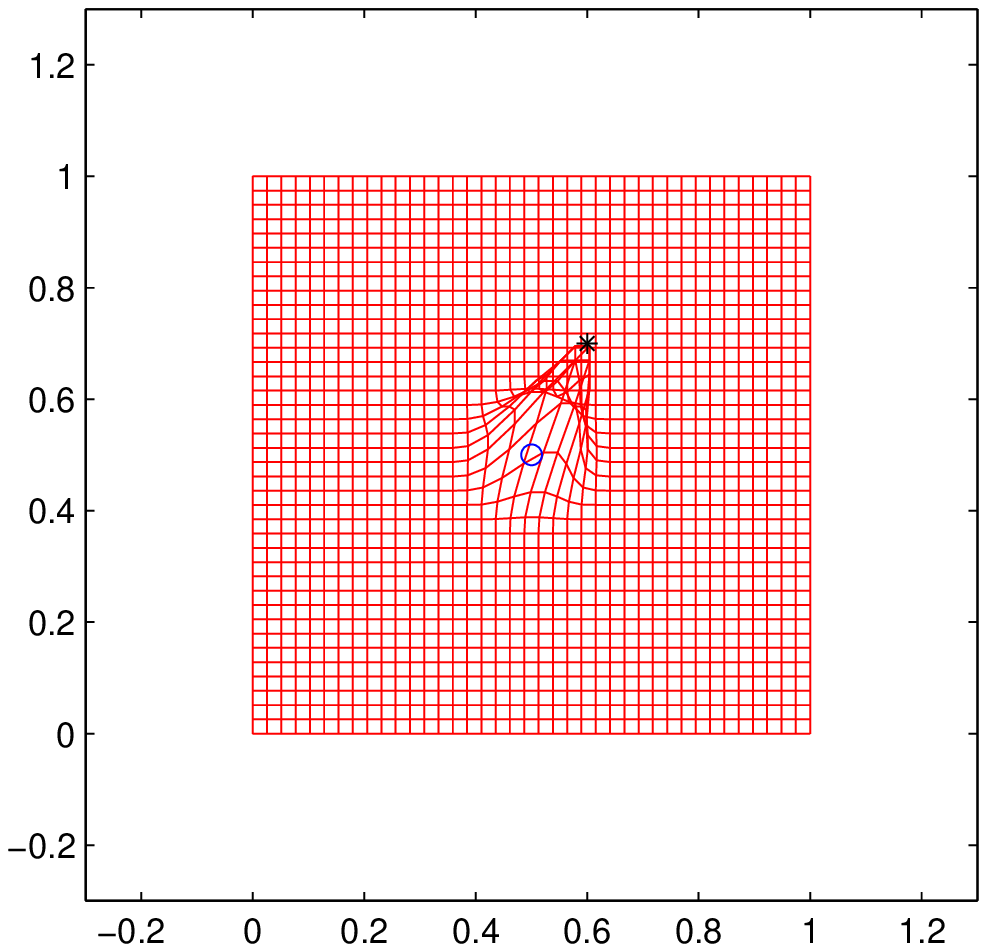}
\centerline{(b) Wu $\psi_{1,2}$, $c=0.15$}
\end{minipage} 
\begin{minipage}{40mm}
\includegraphics[width=4cm]{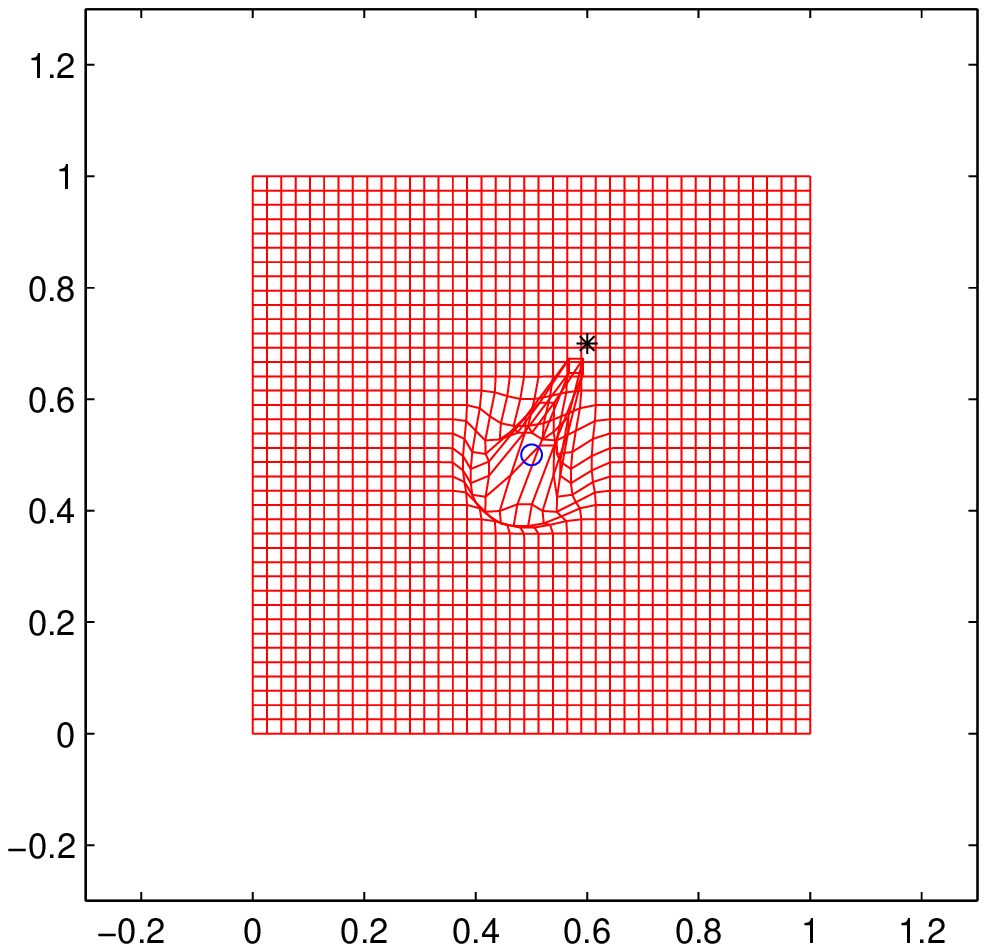}
\centerline{(c) Gneiting $\tau_{2,7/2}$, $c=0.15$}
\end{minipage} 
\begin{minipage}{40mm}
\includegraphics[width=4cm]{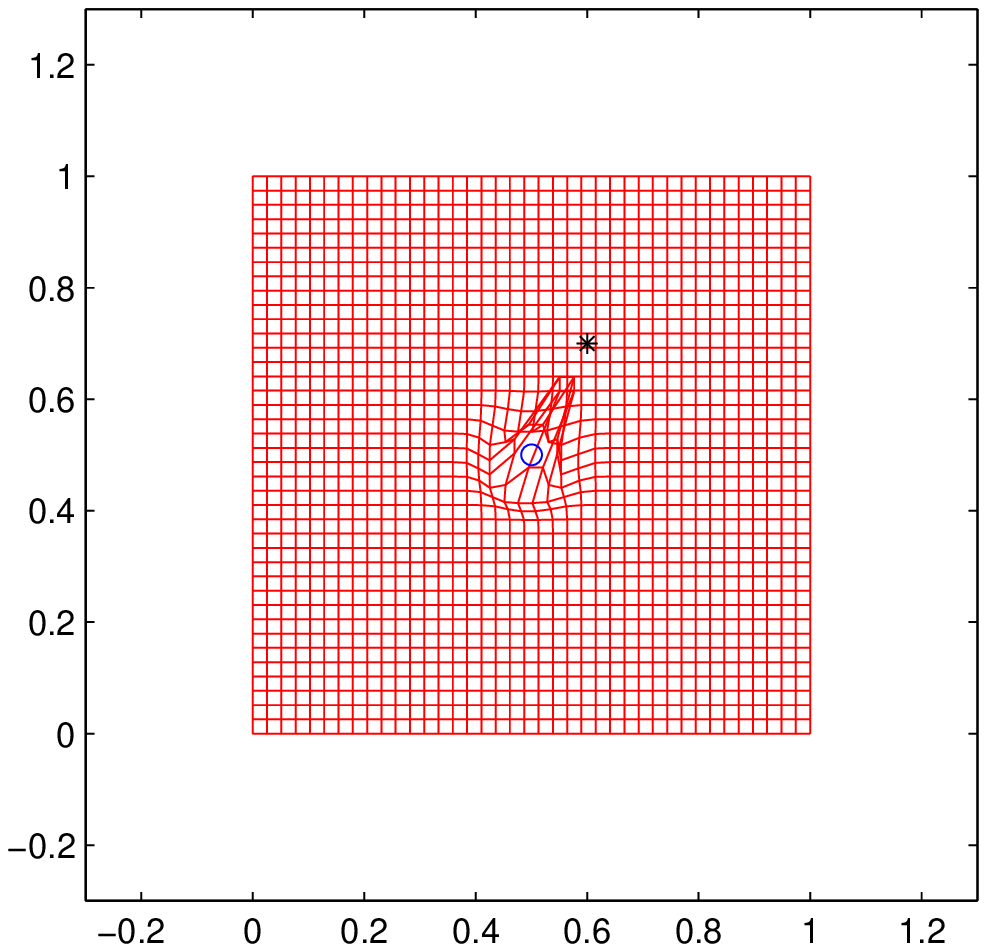}
\centerline{(d) Gneiting $\tau_{2,5}$, $c=0.15$}
\end{minipage}
\caption{Deformation results of one-landmark matching using locality parameters which do not satisfy the topology preservation condition. The source landmark is marked by a circle ($\circ$), while the target one by a star ($\ast$).}
\end{center}
\label{1_landmark_c_no}
\end{figure*}

\section{Topology preservation: the four-landmark case} \label{sec:4}

For more extended deformations, we consider much larger supports which are able to cover whole image. Here, the influence of each landmark extends to the entire domain, therefore global deformations will be generated. In the following we compare topology preservation properties of large extended deformations, and we can set locality parameter large to fulfill this purpose. For this aim, we consider four inner landmarks in a grid, located so as to form a rhombus at the center of the figure, and assuming that only the lower vertex is downward shifted of $\Delta$ \cite{Yang11}. The landmarks of source and target images are
\begin{small}
$$
P=\{P_1=(0,1),\ P_2=(-1,0),\ P_3=(0,-1),\ P_4=(1,0)\}
$$ 
\end{small}
and 
\begin{small}
$$
Q=\{Q_1=(0,1),\ Q_2=(-1,0),\ Q_3=(0,-1-\Delta),\ Q_4=(1,0)\},
$$ 
\end{small}
respectively, with $\Delta>0$.

In this case, we explicitly write the two components of a generic transformation $\textbf{H}:\RR^2\rightarrow \RR^2$ obtained by a transformation of four points $P_1$, $P_2$, $P_3$ and $P_4$, i.e.
\begin{align*}
H_1(\textbf{x})=x+\sum_{i=1}^4 c_{1,i} \Phi(||\textbf{x}-P_i||),
\end{align*}
and
\begin{align*}
H_2(\textbf{x})=y+\sum_{i=1}^4 c_{2,i}\Phi(||\textbf{x}-P_i||).
\end{align*}
We transform $P_i$ to $Q_i$, with $i=1,\dots ,4$, to obtain the coefficients $c_{1,i}$ and $c_{2,i}$. To do that, we need to solve two systems of four equations in four unknowns, whose solutions are given by

\begin{align*}
c_{1,k}=0,\quad k=1,\ldots,4,
\end{align*}
and
\begin{align*}
\begin{array}{rcl}
c_{2,1} & = & \displaystyle{\frac{\beta^2+\beta-2\alpha^2}{(1-\beta)[(1+\beta)^2-4\alpha^2]}\Delta,}\\ 
c_{2,2} & = & \displaystyle{\frac{\alpha}{(1+\beta)^2-4\alpha^2}\Delta,}\\
c_{2,3} & = & \displaystyle{-\frac{1+\beta-2\alpha^2}{(1-\beta)[(1+\beta)^2-4\alpha^2]}\Delta,}\\
c_{2,4} & = &  c_{2,2},
\end{array}
\end{align*}
where 
$$\alpha=\Phi\left(\frac{\sqrt{2}}{c}\right),\hskip 0.2cm \beta=\Phi\left(\frac{2}{c}\right).
$$ 
For simplicity, we denote 
\begin{align*}
\begin{array}{rcl}
\Phi_i&=&\Phi(||(x,y)-P_i||/c),\quad i=1,\ldots,4. 
\end{array}
\end{align*}

Then the determinant of the Jacobian is
\begin{align*}
\det\left(J(x,y)\right)=1+\sum_{i=1}^4 c_{2,i}\frac{\partial\Phi_i}{\partial y}.
\end{align*}
For analyzing the positivity of the determinant of the Jacobian matrix, we calculate the minimum value at position $(0,y)$, with $y>1$. In the following, we analyze the value of the Jacobian determinant at $(0,y)$, with $y>1$, for these four CSRBFs. Because, in this case, we choose a very large parameter $c$, in order to get the value of the Jacobian determinant, we consider $||\cdot||/c$ to be infinitesimal and omit terms of higher order. 

\subsection{Gneiting $\tau_{2,7/2}$} \label{sec:5.1}
The approximation of Gneiting's function $\tau_{2,7/2}$, obtained by using Taylor's expansion, is given by
\begin{align*}
\begin{array}{rcl}
\Phi(r) & = & \displaystyle{(1-r)^{\frac{7}{2}}\left(1+\frac{7}{2}r-\frac{135}{8}r^2\right)} \\
& \approx & \displaystyle{1-\frac{99}{4}r^2+\frac{1155}{16}r^3},
\end{array}
\end{align*}
while its first derivative can be approximated as
\begin{align*}
\Phi'(r)\approx -\frac{99}{2}r+\frac{3465}{16}r^2.
\end{align*}
Thus, the values of $\alpha$ and $\beta$ are
\begin{align*}
\alpha=\Phi\left(\frac{\sqrt{2}}{c}\right)\approx 1-\frac{99}{2}\frac{1}{c^2}+\frac{1155\sqrt{2}}{8}\frac{1}{c^3},
\end{align*}
\begin{align*}
\beta=\Phi\left(\frac{2}{c}\right)\approx 1-\frac{99}{c^2}+\frac{1155}{2}\frac{1}{c^3}.
\end{align*}
Using now numerical estimates of $\alpha$ and $\beta$, we get
\begin{align*}
(1+\beta)^2-4\alpha^2 \approx  1155\left(\frac{2-\sqrt{2}}{c^3}\right),
\end{align*}
from which it follows
\begin{align*}
c_{2,1}\approx-\frac{c^3}{1155(2-\sqrt{2})}\Delta.
\end{align*}
Therefore we can compute the Jacobian determinant at $(0,y)$, for $y>1$, as

\begin{align}\label{det_gau}
\begin{array}{rcl}
\displaystyle{\det\left(J(0,y)\right)} &=& \displaystyle{1+c_{2,1}\left(\left.\frac{\partial\Phi_1}{\partial y}\right|_{(x,y)=(0,y)}-\left.\frac{\partial\Phi_2}{\partial y}\right|_{(x,y)=(0,y)} \right.} \\
&+& \displaystyle{\left.\left.\frac{\partial\Phi_3}{\partial y}\right|_{(x,y)=(0,y)}-\left.\frac{\partial\Phi_4}{\partial y}\right|_{(x,y)=(0,y)}\right),}
\end{array}
\end{align}
evaluating singularly the four partial derivatives 
\begin{align*}
\left.\frac{\partial\Phi_1}{\partial y}\right|_{(x,y)=(0,y)} & \approx \frac{1}{c}\left(-\frac{99}{2}\frac{y-1}{c}+\frac{3465}{16}\frac{(y-1)^2}{c^2}\right), \\
\left.\frac{\partial\Phi_2}{\partial y}\right|_{(x,y)=(0,y)} & \approx  -\frac{y}{c^2}\left(\frac{99}{2}-\frac{3465}{16}\frac{\sqrt{1+y^2}}{c}\right), \\
\left.\frac{\partial\Phi_3}{\partial y}\right|_{(x,y)=(0,y)} & \approx  \frac{1}{c}\left(-\frac{99}{2}\frac{y+1}{c}+\frac{3465}{16}\frac{(y+1)^2}{c^2}\right),
\end{align*}
and
$$\frac{\partial\Phi_4}{\partial y}=\frac{\partial\Phi_2}{\partial y}.$$
Replacing then the four partial derivatives in (\ref{det_gau}), we obtain
\begin{align}\label{gne_7}
\det\left(J(0,y)\right) \approx 1-0.6402\Delta\left(y^2+1-y\sqrt{y^2+1}\right).
\end{align}

\subsection{Gneiting $\tau_{2,5}$} \label{sec:4.2}
Gneiting's function $\tau_{2,5}$ can be approximated as
\begin{align*}
\begin{array}{rcl}
\Phi(r)&=&(1-r)^{5}\left(1+5r-27r^2\right)\\
&\approx& 1-42r^2+175r^3,
\end{array}
\end{align*}
while its first derivative is
\begin{align*}
\Phi'(r)\approx -84r+525r^2.
\end{align*}
Then, $\alpha$ and $\beta$ can approximatively be represented as follows
\begin{align*}
\alpha=\Phi\left(\frac{\sqrt{2}}{c}\right)\approx 1-\frac{84}{c^2}+\frac{350\sqrt{2}}{c^3},
\end{align*}
\begin{align*}
\beta=\Phi\left(\frac{2}{c}\right)\approx 1-\frac{168}{c^2}+\frac{1400}{c^3}.
\end{align*}
In order to evaluate $(1+\beta)^2-4\alpha^2$, from approximations of $\alpha$ e $\beta$ we deduce
\begin{align*}
(1+\beta)^2-4\alpha^2  \approx 2800\left(\frac{2-\sqrt{2}}{c^3}\right),
\end{align*}
so we can obtain
\begin{align*}
c_{2,1}\approx-\frac{c^3}{2800(2-\sqrt{2})}\Delta.
\end{align*}
Referring now to (\ref{det_gau}), we can compute the following four partial derivatives for $\tau_{2,5}$, i.e., 
\begin{align*}
\left.\frac{\partial\Phi_1}{\partial y}\right|_{(x,y)=(0,y)} &\approx \frac{1}{c}\left(-84\frac{y-1}{c}+525\frac{(y-1)^2}{c^2}\right),\\
\left.\frac{\partial\Phi_2}{\partial y}\right|_{(x,y)=(0,y)} &\approx -\frac{y}{c^2}\left(84-525\frac{\sqrt{1+y^2}}{c}\right),\\
\left.\frac{\partial\Phi_3}{\partial y}\right|_{(x,y)=(0,y)} &\approx \frac{1}{c}\left(-84\frac{y+1}{c}+525\frac{(y+1)^2}{c^2}\right),
\end{align*}
and
\begin{align}
\frac{\partial\Phi_4}{\partial y}=\frac{\partial\Phi_2}{\partial y}.
\end{align}
Then, replacing such derivatives in (\ref{det_gau}), we get
\begin{align}\label{gne_5}
\det\left(J(0,y)\right) \approx 1-0.6402\Delta\left(y^2+1-y\sqrt{y^2+1}\right).
\end{align}

\begin{figure}
\centering
\begin{minipage}[t]{0.5\textwidth}
\centering
\includegraphics[width=7.5cm]{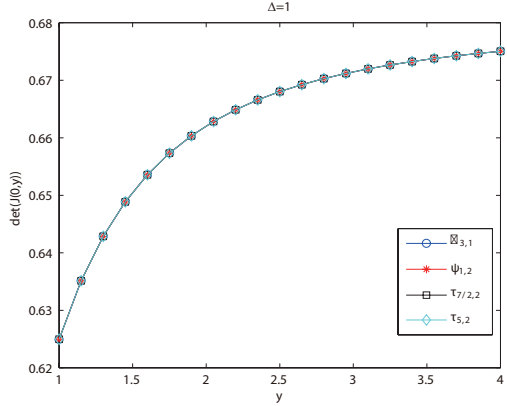}
\end{minipage}
\caption{Value of $\det(J(0,y))$, with $y>1$, by varying CSRBFs.}
\label{global}
\end{figure}

We compare the approximations of (\ref{gne_7}) and (\ref{gne_5}) with the ones acquired by the work \cite{Yang11}. Figure \ref{global} shows the same values of $\det(J(0,y))$, with $y>1$, when one uses as CSRBF transformations based on Wendland's, Wu's and Gneiting's functions. This indicates that functions $\varphi_{3,1}$, $\psi_{2,1}$, $\tau_{2,7/2}$ and $\tau_{2,5}$ have the same behavior. The equations obtained in \cite{Yang11} using Wendland's and Wu's functions and those deduced in (\ref{gne_7}) and (\ref{gne_5}) guarantee the Jacobian determinant is positive for any $y>1$. This means all these transformations can easily preserve topology. 

\subsection{Numerical experiments} \label{sec:5.3}

In this subsection, a schematic diagram of four landmarks and a real case of brain images are evaluated. Firstly, we consider a grid $[0,1]\times[0,1]$ and compare results obtained by its distortion, which is created by the shift of one of the four landmarks distributed in rhomboidal position. The source landmarks are 
$$
\{(0.5,0.65), (0.35,0.5), (0.65,0.5), (0.5,0.35)\}
$$ 
and are respectively transformed in the following target landmarks 
$$
\{(0.5,0.65),(0.35,0.5), (0.65,0.5), (0.5,0.25)\} .
$$
Taking $c=100$ as support size, we obtain Figure 5.2 in case of CSRBFs. 

\begin{figure*} 
\begin{center}
\begin{minipage}{40mm}
\centering
\includegraphics[width=4cm]{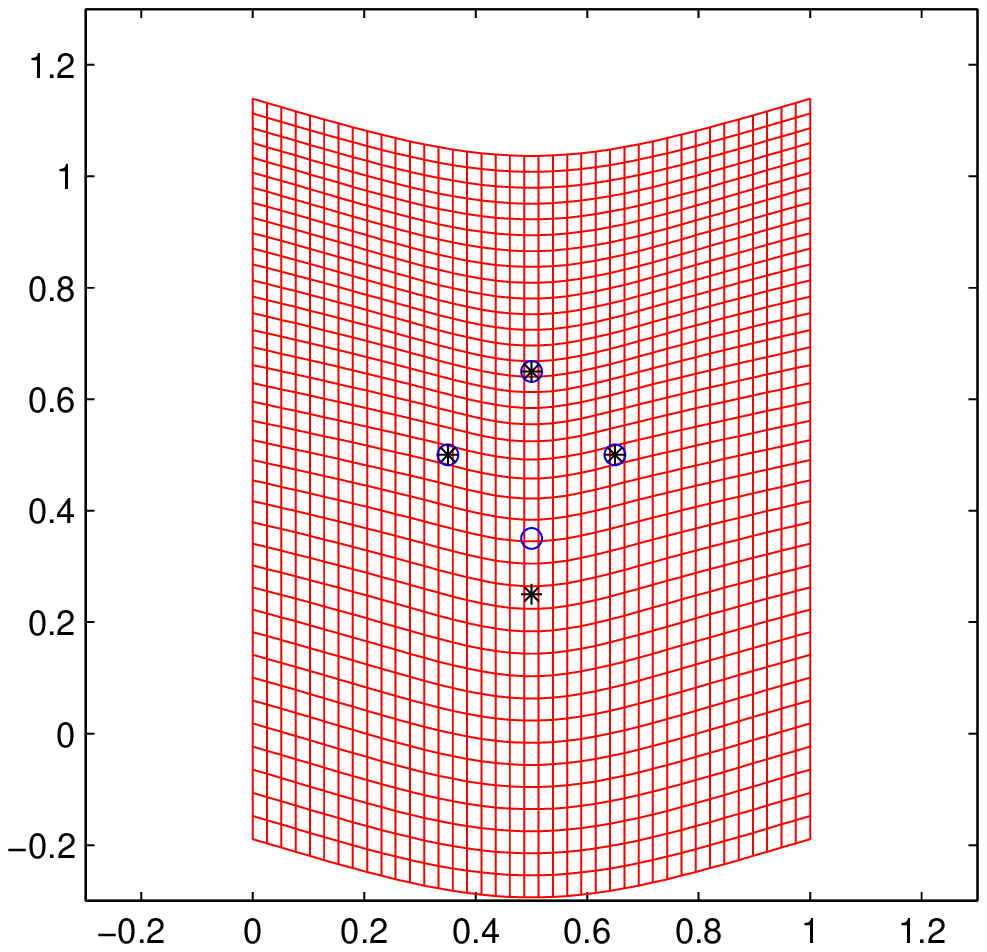}
\centerline{(a) Wendland $\varphi_{3,1}$, $c = 100$}
\end{minipage} 
\begin{minipage}{40mm}
\centering
\includegraphics[width=4cm]{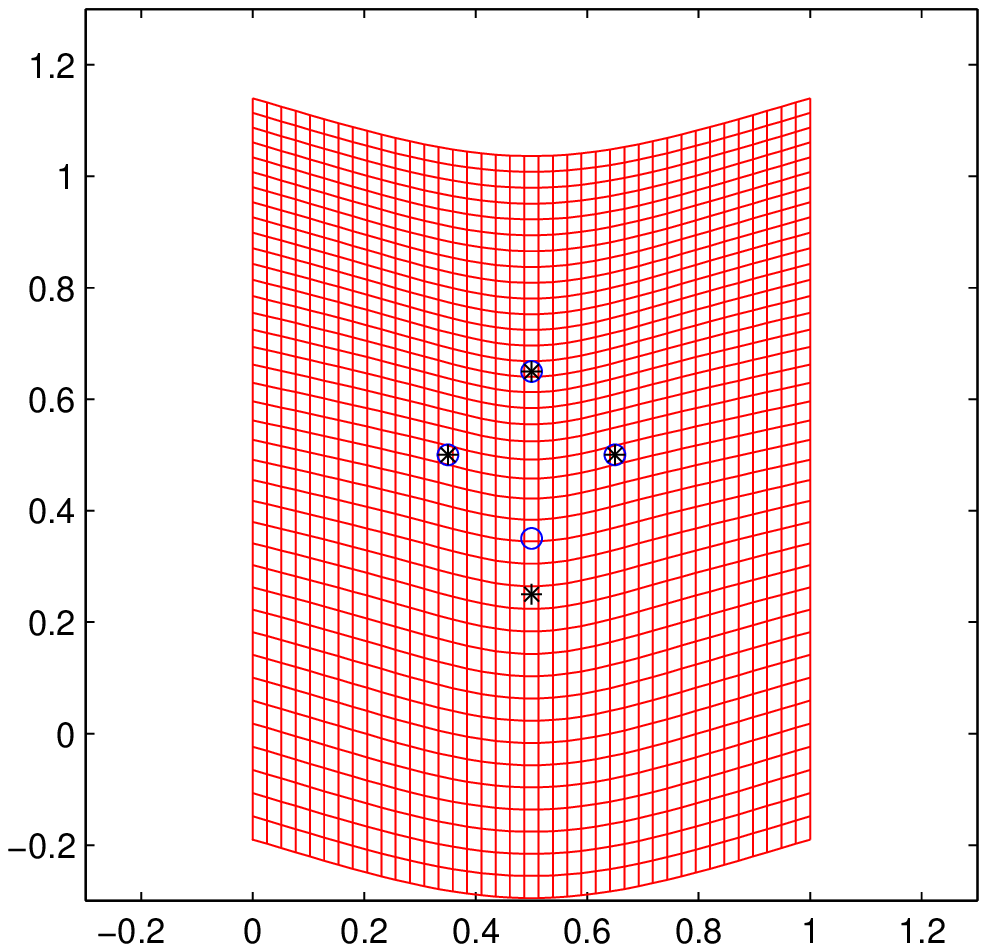}
\centerline{(b) Wu $\psi_{1,2}$, $c = 100$}
\end{minipage} 
\begin{minipage}{40mm}
\centering
\includegraphics[width=4cm]{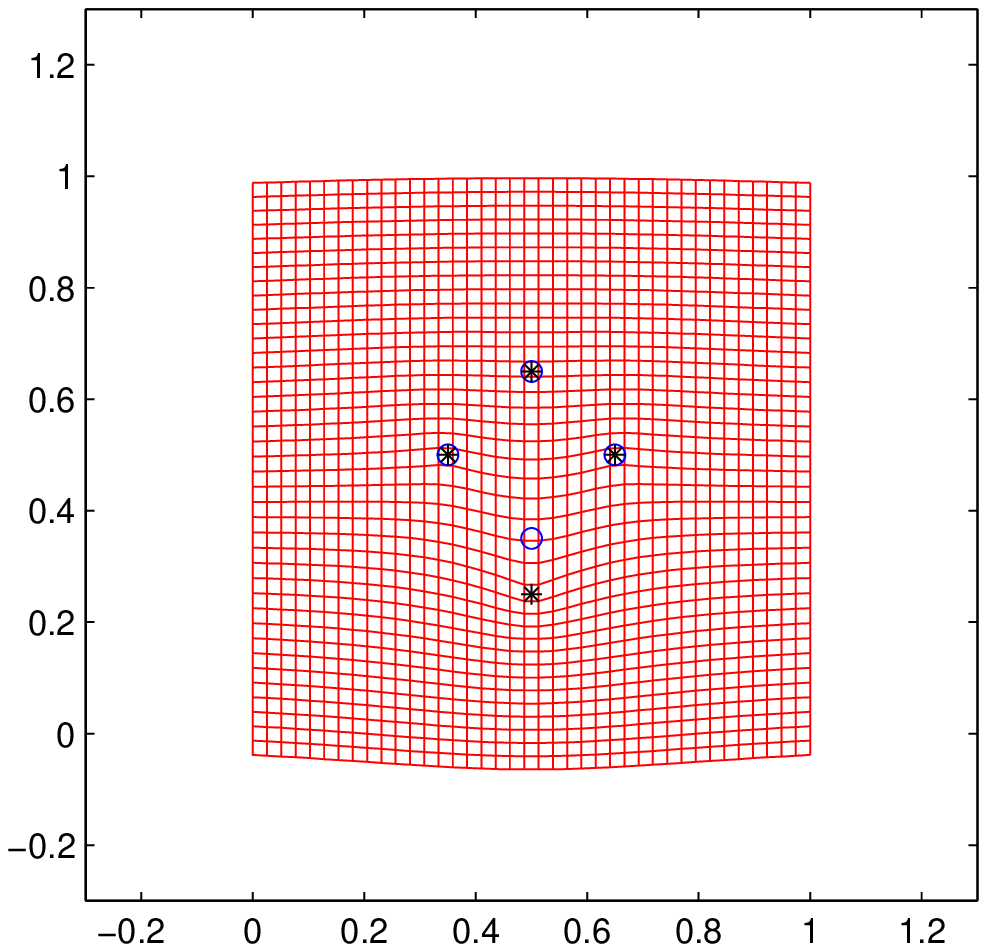}
\centerline{(c) Gneiting $\tau_{2,7/2}$, $c = 100$}
\end{minipage} 
\begin{minipage}{40mm}
\centering
\includegraphics[width=4cm]{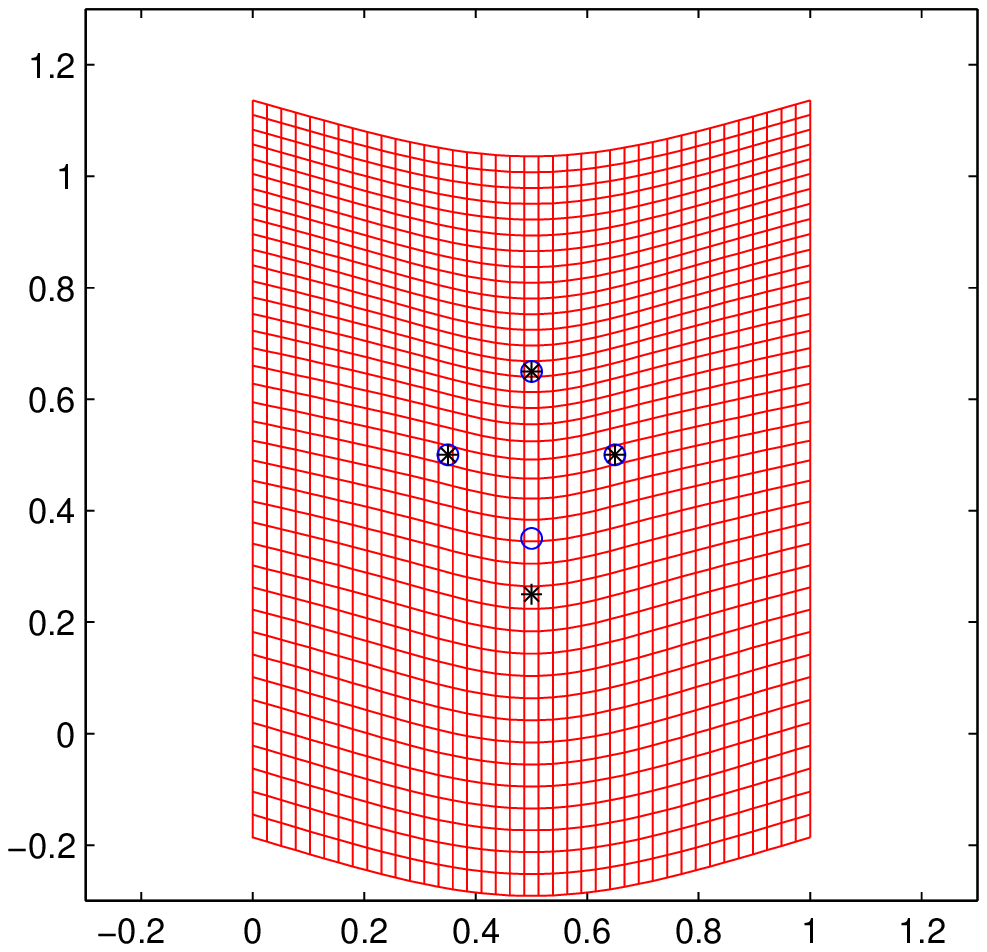}
\centerline{(d) Gneiting $\tau_{2,5}$, $c = 100$}
\end{minipage} 
\caption{Deformation results of four landmarks; the source landmarks are marked by a circle ($\circ$), while the target ones by a star ($\ast$).}
\end{center}
\label{4punti}
\end{figure*}

In agreement with theoretical results, Figure 5.2 confirms that all these functions can preserve topology and all of them present very similar deformations, with the exception of the $\tau_{2,7/2}$ function. In this case, we can see that, even if the support is very large, the image is only deformed slightly on local area. This is the best result obtained. 

We can conclude that among these four functions, $\tau_{2,7/2}$ can lead to good result in case of high landmarks density, i.e. when distance among landmarks is very small and each of them influences the whole image. Furthermore, it not only preserves topology, but also changes shape very slightly in the whole image when support is relatively large. 

Actually, \cite{Cavoretto13} reported the properties of $\tau_{2,7/2}$ and $\tau_{2,5}$ functions in the case of a large number of landmarks for square shift and scaling, and also for circle contraction and expansion. Numerical results showed that these two functions have good performances in those cases. 

In Figure 5.3, we show the case for $c=0.15$ in which topology preservation is not satisfied. The transformed images are deeply misrepresented mainly close to the shifted point.

\begin{figure*} 
\begin{center}
\begin{minipage}{40mm}
\includegraphics[width=4cm]{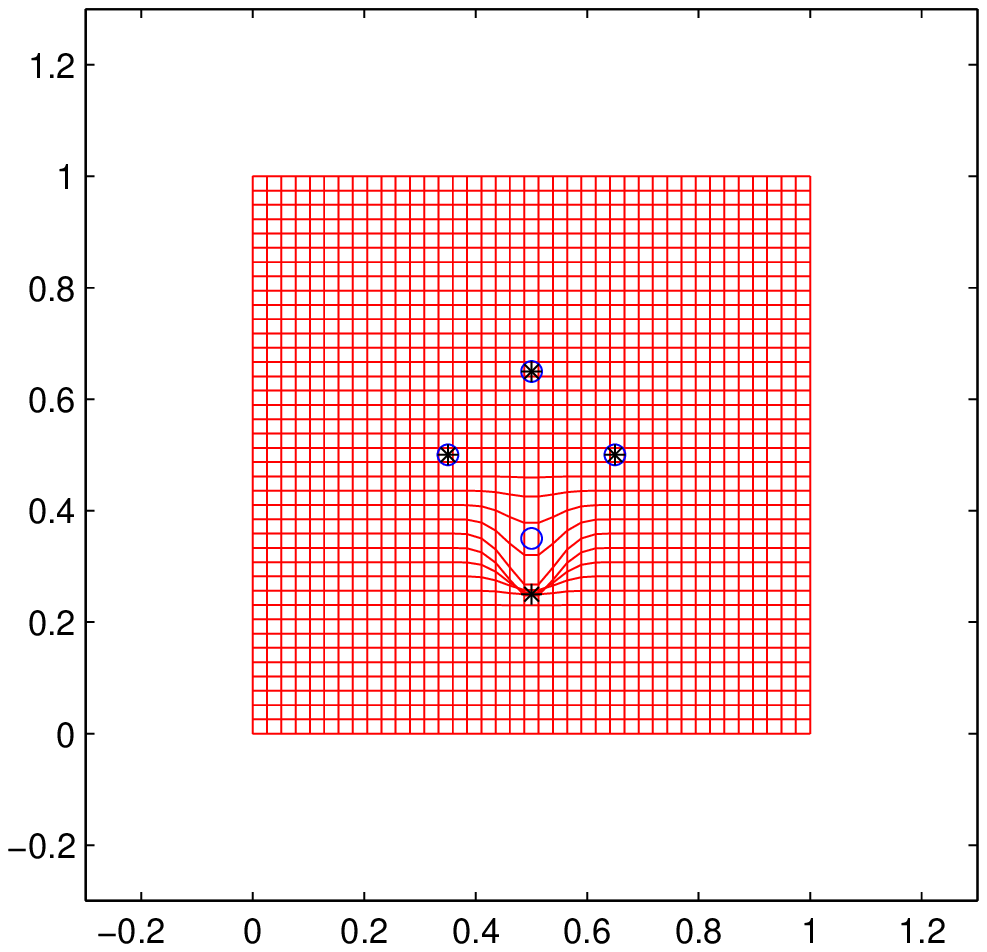}
\centerline{(a) Wendland $\varphi_{3,1}$, $c = 0.15$}
\end{minipage} 
\begin{minipage}{40mm}
\includegraphics[width=4cm]{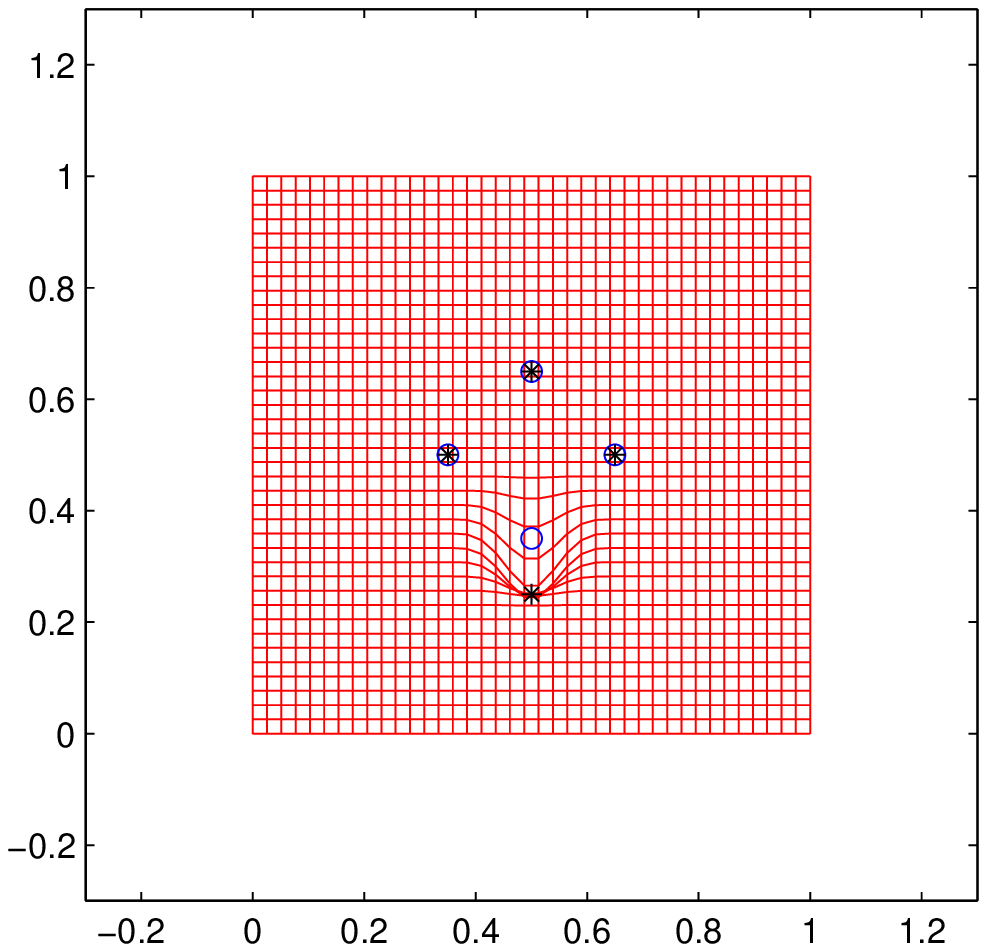}
\centerline{(b) Wu $\psi_{1,2}$, $c = 0.15$}
\end{minipage} 
\begin{minipage}{40mm}
\includegraphics[width=4cm]{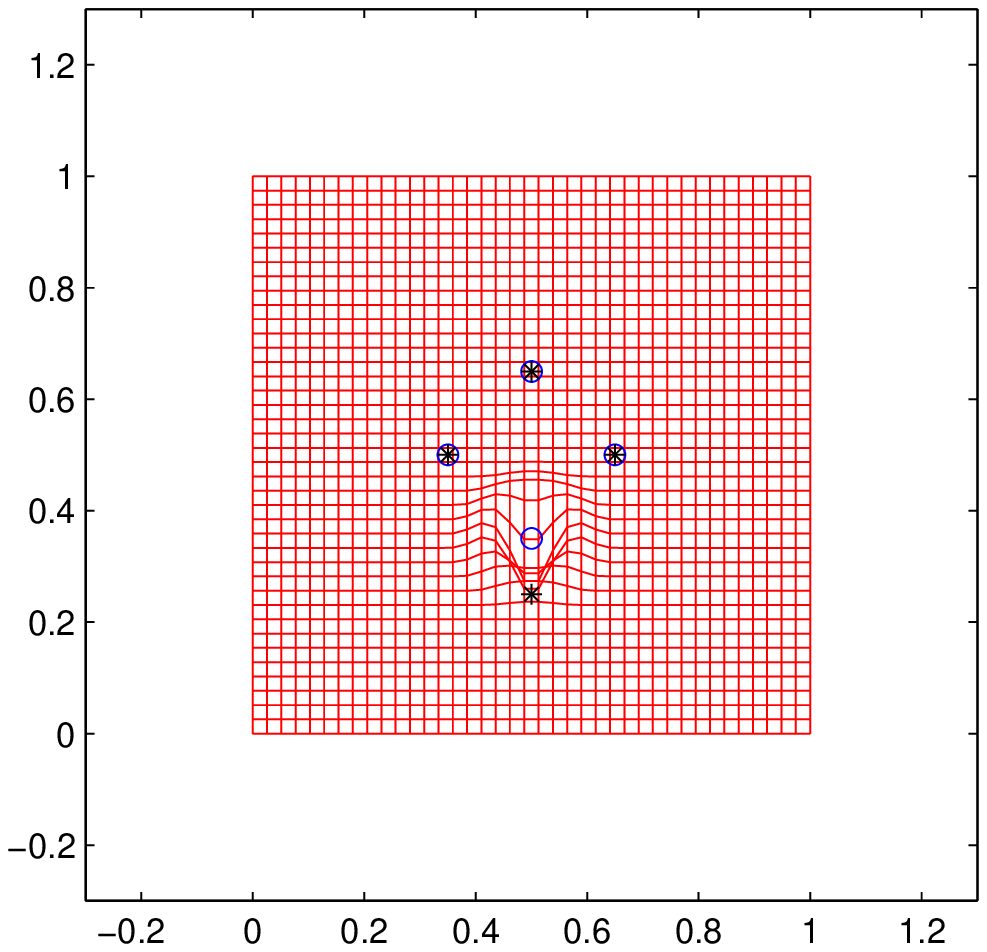}
\centerline{(c) Gneiting $\tau_{2,7/2}$, $c = 0.15$}
\end{minipage} 
\begin{minipage}{40mm}
\includegraphics[width=4cm]{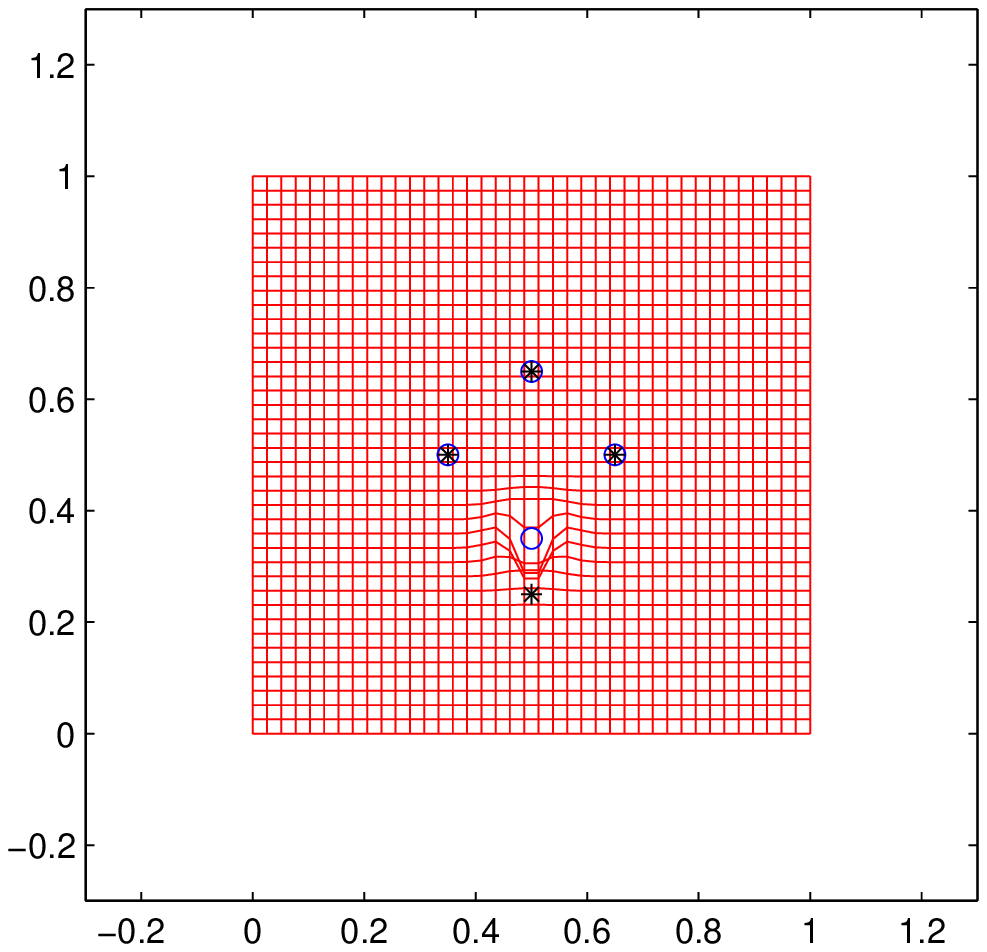}
\centerline{(d) Gneiting $\tau_{2,5}$, $c = 0.15$}
\end{minipage} 
\caption{Deformation results of four landmarks using locality parameters which do not satisfy the topology preservation condition; the source landmarks are marked by a circle ($\circ$), while the target ones by a star ($\ast$).}
\end{center}
\label{4punti_no}
\end{figure*} 

Moreover, a real medical case of landmark-based image registration is displayed in Figure 5.4, which focuses on two brain images. The image (a) is the source image with the corresponding landmarks marked by \textcolor{red}{$\circ$}, whereas (b) is the target image with the landmarks \textcolor{green}{$\ast$}. In this case, we choose two different support sizes for $\tau_{2,7/2}$ in order to observe the various topology behaviours. More precisely, when $c=20,$ which is relatively a large value, the transformed image can preserve topology well. However, the opposite situation occurs when $c$ is  smaller, for instance equal to 2.

\begin{figure*} 
\begin{center}
\begin{minipage}{90mm}
\includegraphics[width=9cm]{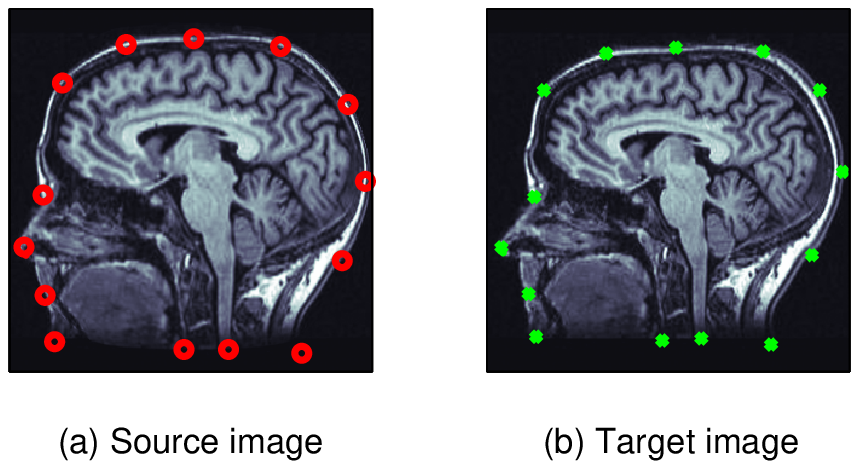}
\end{minipage} \hskip -1cm
\begin{minipage}{90mm}
\includegraphics[width=9cm]{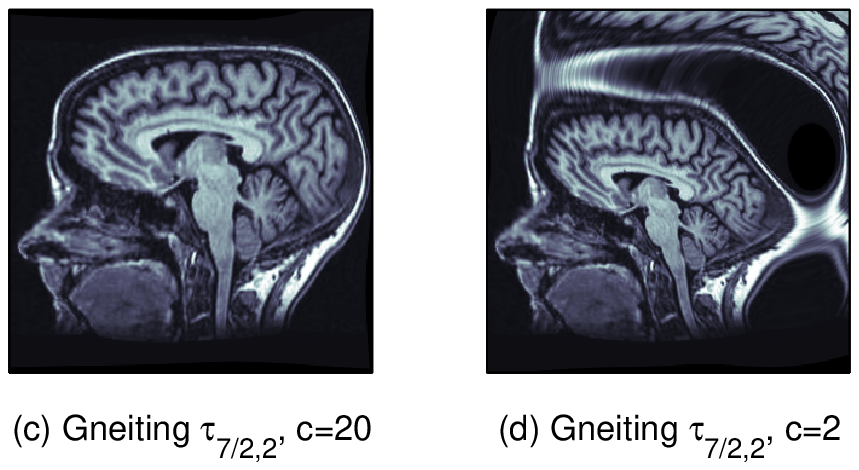}
\end{minipage} \vskip -0.3cm
\caption{Registration results of brain images using transformation $\tau_{2,7/2}$ with various support sizes $c.$ The source landmarks are marked by a circle (\textcolor{red}{$\circ$}), while the target ones by a star (\textcolor{green}{$\ast$}).}
\end{center}
\label{brain_case}
\end{figure*}

\section{Conclusions}

For image registration, we should guarantee that images preserve their original structure and they are not folded after deformation. Hence the transformations we use should be preserved topologically. 

In this paper, we evaluated the topology preservation of two kinds of Gneiting's functions, and compared the results with Wendland's and Wu's functions in one-, four-landmark and medical brain image cases. In the first case, all functions have very similar performances. Conversely, in four-landmark and brain image cases, Gneiting's functions have better performances, especially the $\tau_{2,7/2}$ function provides the best registration result.


\section*{Acknowledgements} 
The authors are very grateful to the anonymous referees for their detailed and valuable comments which helped to greatly improve the paper. The second and third authors acknowledge support from the GNCS-INdAM. 



\end{document}